\documentclass[12pt]{article}
\usepackage{a4wide}
\usepackage{amssymb, amsmath, amsthm}

{
\bigskip
}

\theoremstyle{theorem}
\newtheorem{st}{Theorem}
\newtheorem{lem}{Lemma}
\newtheorem{fa}{Fact}
\theoremstyle{definition}
\newtheorem*{bewijs}{Proof}

\newcommand{\bew}[1]{\begin{bewijs} #1 \qed \end{bewijs}}

\begin{document}

\title{The relative sizes of sumsets and difference sets}
\author{Merlijn Staps}
\date{October 9, 2014}
\maketitle

\begin{abstract}\noindent
Let $A$ be a finite subset of a commutative additive group $Z$. The sumset and difference set of $A$ are defined as the sets of pairwise sums and differences of elements of $A$, respectively. The well-known inequality $\sigma(A)^{1/2} \leq \delta(A) \leq \sigma(A)^2,$ where $\sigma(A)=\frac{|A+A|}{|A|}$ is the doubling constant of A and $\delta(A)=\frac{|A-A|}{|A|}$ is the difference constant of A, relates the relative sizes of the sumset and difference set of A. The exponent $2$ in this inequality is known to be optimal, for the exponent $\frac12$ this is unknown. We determine those sets for which equality holds in the above inequality. We find that equality holds if and only if $A$ is a coset of some finite subgroup of $Z$ or, equivalently, if and only if both the doubling constant and difference constant are equal to $1$. This implies that there is space for possible improvement of the exponent $\frac12$ in the inequality. We then use the derived methods to show that Pl\"unnecke's inequality is strict when the doubling constant is larger than $1$.\end{abstract}

Let $(Z,+)$ be a commutative group. We consider finite non-empty subsets $A \subset Z$. The sumset $A+A$ and difference set $A-A$ of $A$ are defined by $A+A = \{a+b : a,b \in A\}$ and $A-A = \{a-b: a,b \in A\}$. More generally we define
\[
nA-mA = \{a_1 + \cdots + a_n - b_1 - \cdots - b_m : a_1, \ldots, a_n, b_1, \ldots, b_m \in A\}.
\]
Furthermore, we define the the doubling constant $\sigma(A) = \frac{|A+A|}{|A|}$ and difference constant $\delta(A) = \frac{|A-A|}{|A|}$ of $A$. It is clear that $\sigma(A) \geq 1$ and $\delta(A) \geq 1$. Of interest are those sets $A$ for which $\sigma(A)$ or $\delta(A)$ are ``small''. When $\sigma(A)=1$ or $\delta(A)=1$ we have the following result that is easily proven (for instance, see \cite{TV}):

\begin{fa}\label{fa:sigma1}
For a set $A$ we have $\sigma(A) = 1$ if and only if $\delta(A)=1$ if and only if $A$ is a coset of some finite subgroup of $Z$.
\end{fa}

It is a general result that sets with a small sumset also have small difference set and vice versa. In particular, it turns out that the following inequality
\begin{align}
\sigma(A)^{1/2} \leq \delta(A) \leq \sigma(A)^2 \label{eq:ineq}
\end{align}
relates the doubling constant $\sigma(A)$ and difference constant $\delta(A)$ \cite{TV}.

The bounds in \eqref{eq:ineq} are the best known bounds of this type. The exponent $2$ in the upper bound cannot be improved at all \cite{Hennecart}. Whether the exponent $\frac12$ in the lower bound can be improved is not known. Here, we determine the equality case for both inequalities in \eqref{eq:ineq}. The main result of the paper is the following theorem.

\begin{st}\label{th:main-theorem}
We have $\sigma(A) = \delta(A)^2$ or $\delta(A)=\sigma(A)^2$ if and only if $A$ is a coset of a finite subgroup of $Z$, i.e. if and only if $\sigma(A)=\delta(A)=1$.
\end{st}

In section \ref{sec:upper} we prove the known inequality $\delta(A) \leq \sigma(A)^2$, the upper bound in \eqref{eq:ineq}, which easily follows from Ruzsa's triangle inequality \cite{Ruzsa2}. We show that equality holds if and only if $\sigma(A)=\delta(A)=1$ (theorem \ref{th:upper}).

In section \ref{sec:lower} we determine the equality case of the lower bound in \eqref{eq:ineq}. We first derive an equality condition for a lemma of Petridis \cite{Petridis} that can be used to prove the lower bound. We then determine that equality holds in the lower bound of \eqref{eq:ineq} if and only if $\sigma(A)=\delta(A)=1$ (theorem \ref{th:lower}). The fact that equality holds only in this case is a necessary condition for a possible improvement of the exponent $\frac12$ in \eqref{eq:ineq}.



Petridis' lemma can also be used to derive Pl\"unnecke's inequality \cite{Petridis}. This inequality states $|nA| \leq \sigma(A)^n|A|$ and thus gives an upper bound on the size of sumsets of the form $A+A+\cdots+A$ \cite{Plunnecke}. In section \ref{sec:Plunnecke} we use the results derived in section \ref{sec:lower} to show that Pl\"unnecke's inequality is strict unless $\sigma(A)=1$.

We will use the symbols $\subset$, $\subsetneq$ and $\sqcup$ for inclusion, strict (proper) inclusion and disjoint union, respectively.

\section{Sets with few sums and many differences\label{sec:upper}}

We consider the inequality $\delta(A) \leq \sigma(A)^2$. Careful analysis of the standard proof using the Ruzsa triangle inequality \cite{Ruzsa2} shows that the equality case of this inequality is given by those sets $A$ for which $\delta(A)=\sigma(A)=1$.

\begin{st}\label{th:upper}
For $A \subset Z$ we have $\delta(A) \leq \sigma(A)^2$ with equality if and only if $\sigma(A)=1$.
\end{st}

\bew{The inequality $\delta(A) \leq \sigma(A)^2$ is a special case of a more general inequality
\[
|K||J-L| \leq |J-K||K-L|
\]
which is known as the Ruzsa triangle inequality \cite{Ruzsa2}. The inequality is proven by constructing an injective map $K \times (J-L) \to (J-K) \times (K-L)$. We construct this map in the special case $(J,K,L) = (A,-A,A)$, thereby proving $|A||A-A| \leq |A+A|^2$. First we choose a map $\psi: A-A \to A^2$ such that we have $u-v=w$ when $\psi(w)=(u,v)$. Now consider the map $\phi: A \times (A-A) \to (A+A)^2$ such that $\phi(a,u) = (a+b,a+c)$ where $(b,c) = \psi(u)$. This map $\phi$ is easily proven to be injective.

We have equality in $\delta(A)\leq \sigma(A)^2$ if and only if the map $\phi$ is also surjective. Suppose $\phi$ is surjective. Then for each $k \in A+A$ there exists a pair $(a,u) \in A \times (A-A)$ with $\phi(a,u) = (k,k)$. It follows that $u=0$ and $a = k - b$ where $b$ is the first coordinate of $\psi(0)$. We conclude that $A+A$ is contained in $b+A$, hence $|A+A| \leq |b+A| = |A|$ and $\sigma(A) \leq 1$. It follows that $\sigma(A)=1$. Since for $\sigma(A)=1$ we also have $\sigma(A)^2 = 1 = \delta(A)$, the proof is complete.}

\section{Sets with few differences and many sums\label{sec:lower}}

We first state and prove the lemma that is used to prove the lower bound in \eqref{eq:ineq}.

\begin{lem}[Petridis, \cite{Petridis}]\label{le:Petridis}
Let $A, X \subset Z$ such that $|A+X|=K|X|$ and $|A+X'| \geq K|X'|$ for all subsets $X'$ of $X$. Then we have $|A+X+C| \leq K|X+C|$ for all $C \subset Z$.
\end{lem}

\bew{Write $C = \{c_1,\ldots,c_m\}$ and let $C_k = \{c_1,\ldots,c_k\}$. We show $|A+X+C_k| \leq K|X+C_k|$ by induction, the base case $k=1$ being trivial. Notice that
\begin{align}
X+A+C_k = (X+A+C_{k-1}) \cup ((X+A+c_k) \backslash (X_k+A+c_k)) \label{eq:union}
\end{align}
where $X_k = \{x \in X: x+A+c_k \subset X+A+C_{k-1}\}$. Since $X_k +A +c_k \subset X+A+C_k$ we have
\begin{align*}
|X+A+C_k| &\leq |X+A+C_{k-1}| + |X+A+c_k| - |X_k + A+ c_k| \\ &= |X+A+C_{k-1}| + |X+A| - |X_k+A|,
\end{align*}
with equality if and only if the union in \eqref{eq:union} is disjoint. Using the induction hypothesis we now find $|X+A+C_k| \leq K(|X+C_{k-1}| + |X| - |X_k|)$. We also used that $|X_k+A| \geq K|X_k|$, which follows from the fact that $X_k$ is a subset of $X$. Notice that $X+C_k = (X+C_{k-1}) \sqcup ((X+c_k) \backslash (Y_k + c_k))$ where $Y_k = \{x \in X: x+c_k \in X+C_{k-1}\}$. It follows that $|X+C_k| = |X+C_{k-1}| + |X| - |Y_k| \geq |X+C_{k-1}| + |X| - |X_k|$ as $Y_k \subset X_k$, hence
\[
|X+A+C_k| \leq K(|X+C_{k-1}| + |X| - |X_k|) \leq K|X+C_k|,
\]
completing the induction step.
}

We will call two sets $A$ and $B$ \emph{independent} when all sums $a+b$ are different for $a \in A$, $b \in B$, i.e. when $|A+B|=|A||B|$. Using this definition, we can formulate conditions for a set $C$ to satisfy $|A+X+C|=K|X+C|$ in the above lemma. It turns out that such a set $C$ contains a set $Q$ such that $A+X$ and $Q$ are independent and $X+C=X+Q$. This means that the elements of $C$ that are in $Q$ introduce only new elements on both sides of $|A+X+C|=K|X+C|$ whereas the elements that are not in $Q$ introduce no new elements. The subset $Q$ is in general not unique.

\begin{lem}[Equality case of lemma \ref{le:Petridis}]\label{le:Petridis-eq}
Let $A, X \subset Z$ such that $|A+X|=K|X|$ and $|A+X'|>K|X'|$ for all proper non-empty subsets $X'$ of $X$. Then we have $|A+X+C| \leq K|X+C|$ for all $C \subset Z$, with equality if and only if there exists a subset $Q \subset C$ such that $X+C = X+Q$ and such that $A+X$ and $Q$ are independent.
\end{lem}

\bew{First suppose $X+C=X+Q$ and $A+X$ and $Q$ are independent. This implies that $X$ and $Q$ are independent as well and that $A+X+C=A+X+Q$. We now have $|A+X+C| = |A+X+Q| = |A+X||Q| = K|X||Q| = K|X+Q| = K|X+C|$.

Now suppose $|A+X+C|=K|X+C|$. Then we have equality in lemma \ref{le:Petridis}, thus for each $1 \leq k \leq m$ the following conditions are satisfied:
\begin{itemize}
\item we have $(X+A+C_{k-1}) \cap ((X+A+c_k) \backslash (X_k+A+c_k)) = \emptyset$ (since the union in \eqref{eq:union} is disjoint);
\item we have $X_k = \emptyset$ or $X_k = X$ (since we have equality in $|X_k+A| \geq K|X_k|$);
\item we have $Y_k = X_k$ (since we have equality in $|Y_k| \leq |X_k|$).
\end{itemize}
When $X_k = \emptyset$ it follows from the first condition that $(X+A+C_{k-1}) \cap (X+A+c_k) = \emptyset$, hence $X+A+C_k = (X+A+C_{k-1}) \sqcup (X+A+c_k)$ and $(X+C_k) = (X+C_{k-1}) \sqcup (X+c_k)$. When $X_k = X$ we have $Y_k = X$, hence $X+c_k \subset X+C_{k-1}$. It now follows that $X+C_k = X+C_{k-1}$ and $A+X+C_k = A+X+C_{k-1}$. Let $Q$ be the subset of $C$ consisting of those $c_k$ for which $X_k = \emptyset$. Then we have $X+C = \bigsqcup_{q \in Q} (X+q)$, thus $X+C=X+Q$. Furthermore, we have $$|A+X+Q| = |A+X+C| = \left|\bigsqcup_{q \in Q} (A+X+q)\right| = |A+X||Q|,$$ showing that $A+X$ and $Q$ are independent.}

We are now ready to prove the inequality $\sigma(A) \leq \delta(A)^2$ and to determine the equality case.

\begin{st}\label{th:lower}
For $A \subset Z$ we have $\sigma(A) \leq \delta(A)^2$ with equality if and only if $\sigma(A)=1$.
\end{st}

\bew{
Choose the smallest possible non-empty subset $X \subset -A$ minimizing $\frac{|A+X|}{|X|}$. Denote $K = \frac{|A+X|}{|X|} \leq \frac{|A-A|}{|-A|} = \delta(A)$. Then the condition of lemma \ref{le:Petridis-eq} is satisfied, hence for $C=A$ we have $|2A| \leq |2A+X| \leq K|X+A| = K^2|X| \leq K^2|A| \leq \delta(A)^2|A|$, showing that $\sigma(A) \leq \delta(A)^2$. We have equality if there exists a subset $Q \subset A$ such that $X+A=X+Q$ and such that $A+X$ and $Q$ are independent. In that case, it follows that $\delta(A)|X+A| = K|X+A| = |2A+X| = |A+X+Q| = |A+X||Q|$ hence $|Q| = \delta(A)$. Since $A+X$ and $Q$ are independent the sets $A$ and $Q$ are independent, which implies that $|Q|=1$. Thus we have $\delta(A) = |Q| = 1$, implying $\sigma(A)=1$ as well (fact \ref{fa:sigma1}). When $\sigma(A)=1$ we have $\sigma(A) = 1 = \delta(A)^2$.
}

It remains unknown whether the inequality $\sigma(A) \leq \delta(A)^C$ is true for some exponent $C<2$. A necessary condition for this is that the equality in $\sigma(A) \leq \delta(A)^2$ only holds when $\sigma(A)=\delta(A)=1$, which we showed (Theorem \ref{th:lower}).

Penman and Wells \cite{Penman} have shown that $C$ cannot be decreased below $\frac{\log(32/5)}{\log(26/5)} = 1.12594$. There are no better lower bounds for $C$ known.


\section{Pl\"unnecke's inequality is strict when $\sigma(A)>1$\label{sec:Plunnecke}}

Lemma \ref{le:Petridis} can be used to prove Pl\"unnecke's inequality \cite{Petridis,Sanders}. Here we use lemma \ref{le:Petridis-eq} to show that Pl\"unnecke's inequality is strict except when $\sigma(A)=1$.

\begin{st}[Strict Pl\"unnecke inequality]
Suppose that $\sigma(A)>1$. Then we have $|nA| < \sigma(A)^n|A|$ for all $n \geq 1$.
\end{st}

\bew{
The statement is trivially true for $n=1$, so suppose that $n \geq 2$.
Choose the smallest possible non-empty subset $X \subset A$ minimizing $\frac{|A+X|}{|X|}$. Then we have $K = \frac{|A+X|}{|X|} \leq \frac{|A+A|}{|A|} \leq \sigma(A)$. Applying lemma \ref{le:Petridis-eq} with $C=(n-1)A$, we find that
\[
|nA+X| \leq K|(n-1)A + X|.
\]
By induction it follows that $|nA+X|\leq K^n|X|$, yielding
\[
|nA| \leq |nA+X| \leq K^n|X| \leq K^n |A| \leq \sigma(A)^n|A|.
\]
We show that we cannot have equality. If we would have equality, we would have equality in $|2A+X| \leq K|A+X|$, which implies the existence of a $Q \subset A$ such that $|2A+X|=|A+X+Q|=|A+X||Q|$ hence $|Q|=K$. Furthermore, we have $|Q|=1$ since $A+X$ and $Q$ are independent. \\ It follows that $\sigma(A)=K=1$, contradicting the assumption $\sigma(A)>1$.
}

In other words, we have equality in Pl\"unnecke's inequality if and only if $A$ is a coset of some finite subgroup of $Z$.


\vspace{1cm}

Mathematisch Instituut, Universiteit Utrecht, Postbus 80.010, 3508 TA Utrecht, Nederland \\
\emph{E-mail address}: \verb"m.staps@uu.nl"

\end{document}